\documentclass[11pt, reqno, letterpaper]{amsart}
\usepackage{amsfonts}
\usepackage{amsmath}
\usepackage{amssymb}
\usepackage{amscd}
\usepackage[mathscr]{eucal}
\usepackage{url}

\allowdisplaybreaks[1]

\newcommand{\ph}[2]{{\left({#1}\right)}_{#2}}

\newcommand{\gfp}[1]{\Gamma_p{\left({#1}\right)}}
\newcommand{\biggfp}[1]{\Gamma_p{\bigl({#1}\bigr)}}

\newcommand{\biggbin}[2]{\biggl({\genfrac{}{}{0pt}{}{#1}{#2}}\biggr)}

\theoremstyle{plain}
\newtheorem{theorem}{Theorem}[section]

\theoremstyle{definition}

\numberwithin{equation}{section}

\begin{document}

\title[Binomial Coefficient--Harmonic Sum Identities associated to Supercongruences]{Binomial Coefficient--Harmonic Sum Identities\\ associated to Supercongruences}
\author{Dermot M\lowercase{c}Carthy}  

\address{Department of Mathematics, Texas A\&M University, College Station, TX 77843-3368, USA}

\email{mccarthy@math.tamu.edu}


\subjclass[2010]{Primary: 05A19; Secondary: 11B65}


\begin{abstract}
We establish two binomial coefficient--generalized harmonic sum identities using the partial fraction decomposition method. 
These identities are a key ingredient in the proofs of numerous supercongruences. 
In particular, in other works of the author, they are used to establish modulo $p^k$ ($k>1$) congruences between truncated generalized hypergeometric series, and a function which extends Greene's hypergeometric function over finite fields to the $p$-adic setting. 
A specialization of one of these congruences is used to prove an outstanding conjecture of Rodriguez-Villegas which relates a truncated generalized hypergeometric series to the $p$-th Fourier coefficient of a particular modular form.
\end{abstract}

\maketitle

\section{Introduction and Statement of Results}
\makeatletter{\renewcommand*{\@makefnmark}{}
\footnotetext{This work was supported by the UCD Ad Astra Research Scholarship program.}
For non-negative integers $i$ and $n$, we define the generalized harmonic sum, ${H}^{(i)}_{n}$, by
\begin{equation*}
{H}^{(i)}_{n}:= \sum^{n}_{j=1} \frac{1}{j^i}
\end{equation*}
and ${H}^{(i)}_{0}:=0$.
In \cite{C} Chu proves the following binomial coefficient-generalized harmonic sum identity using the partial fraction decomposition method.
If $n$ is a positive integer, then
\begin{equation}\label{AEOZ}
\sum_{k=1}^{n} {\biggbin{n+k}{k}}^2 {\biggbin{n}{k}}^2
 \biggl[ 1+ 2k H_{n+k}^{(1)} + 2kH_{n-k}^{(1)} -4kH_k^{(1)} \biggr]=0.
\end{equation}
This identity had previously been established using the WZ method \cite{AEOZ} and was used by Ahlgren and Ono in proving the Ap\'ery number supercongruence \cite{AO}.
 
In \cite{DMC}, \cite{DMC1} the author establishes various supercongruences between truncated generalized hypergeometric series, and a function which extends Greene's hypergeometric function over finite fields to the $p$-adic setting. 
Specifically, let $p$ be an odd prime and let $n \in \mathbb{Z}^{+}$. For $1 \leq i \leq n+1$, let $\frac{m_i}{d_i} \in \mathbb{Q} \cap \mathbb{Z}_p$ such that $0<\frac{m_i}{d_i}<1$. 
Let $\gfp{\cdot}$ denote Morita's $p$-adic gamma function, $\left\lfloor x \right\rfloor$ denote the greatest integer less than or equal to $x$ and $\langle x \rangle$ denote the fractional part of $x$, i.e. $x- \left\lfloor x \right\rfloor$.
Then define
 \begin{multline*}
{_{n+1}G} \left( \tfrac{m_1}{d_1}, \tfrac{m_2}{d_2}, \dotsc, \tfrac{m_{n+1}}{d_{n+1}} \right)_p\\
:= \frac{-1}{p-1}  \sum_{j=0}^{p-2} 
{\left((-1)^j\biggfp{\tfrac{j}{p-1}}\right)}^{n+1} 
\prod_{i=1}^{n+1} \frac{\biggfp{\langle \frac{m_i}{d_i}-\frac{j}{p-1}\rangle}}{\biggfp{\frac{m_i}{d_i}}}
(-p)^{-\lfloor{\frac{m_i}{d_i}-\frac{j}{p-1}}\rfloor}.
\end{multline*}
Note that when $p \equiv 1 \pmod {d_i}$ this function recovers Greene's hypergeometric function over finite fields. For a complex number $a$ and a non-negative integer $n$ let $\ph{a}{n}$ denote the rising factorial defined by
\begin{equation*}\label{RisFact}
\ph{a}{0}:=1 \quad \textup{and} \quad \ph{a}{n} := a(a+1)(a+2)\dotsm(a+n-1) \textup{ for } n>0.
\end{equation*}
Then, for complex numbers $a_i$, $b_j$ and $z$, with none of the $b_j$ being negative integers or zero, we define the truncated generalized hypergeometric series
\begin{equation*}
{{_rF_s} \left[ \begin{array}{ccccc} a_1, & a_2, & a_3, & \dotsc, & a_r \vspace{.05in}\\
\phantom{a_1} & b_1, & b_2, & \dotsc, & b_s \end{array}
\Big| \; z \right]}_{m}
:=\sum^{m}_{n=0}
\frac{\ph{a_1}{n} \ph{a_2}{n} \ph{a_3}{n} \dotsm \ph{a_r}{n}}
{\ph{b_1}{n} \ph{b_2}{n} \dotsm \ph{b_s}{n}}
\; \frac{z^n}{{n!}}.
\end{equation*}
An example of one the supercongruence results from \cite{DMC1} is the following theorem.
\begin{theorem}[\cite{DMC1} Thm. 2.7] \label{thm_4G2}
Let $r, d \in \mathbb{Z}$ such that $2 \leq r \leq d-2$ and $\gcd(r,d)=1$. Let $p$ be an odd prime such that $p\equiv \pm1 \pmod d$ or $p\equiv \pm r \pmod d$ with $r^2 \equiv \pm1 \pmod d$. If $s(p) := \gfp{\tfrac{1}{d}} \gfp{\tfrac{r}{d}}\gfp{\tfrac{d-r}{d}}\gfp{\tfrac{d-1}{d}}$, then
\begin{align*}
{_{4}G} \Bigl(\tfrac{1}{d} , \tfrac{r}{d}, 1-\tfrac{r}{d} , 1-\tfrac{1}{d}\Bigr)_p
&\equiv
{_{4}F_3} \Biggl[ \begin{array}{cccc} \frac{1}{d}, & \frac{r}{d}, & 1-\frac{r}{d}, & 1-\frac{1}{d}\vspace{.05in}\\
\phantom{\frac{1}{d_1}} & 1, & 1, & 1 \end{array}
\bigg| \; 1 \Biggr]_{p-1}
+s(p)\hspace{1pt} p
\pmod {p^3}.
\end{align*}
\end{theorem}
\noindent A specialization of this congruence is used to prove an outstanding supercongruence conjecture of Rodriguez-Villegas, which relates a truncated generalized hypergeometric series to the $p$-th Fourier coefficient of a particular modular form \cite{DMC},\cite{DMC2}. Similar results to Theorem \ref{thm_4G2} exist for $_4G$ with other parameters, and also $_2G$ and $_3G$.

The main results of the current paper, Theorems \ref{Cor_BinHarId1} and \ref{Cor_BinHarId2} below, are two binomial coefficient--generalized harmonic sum identities which factor heavily into the proofs of all the $_4G$ congruences. Taking particular values for $n, m, l, c_1$ and $c_2$ in these identities allows the vanishing of certain terms in the proofs. Note that letting $m=n$ in Theorem \ref{Cor_BinHarId1} recovers (\ref{AEOZ}).

\begin{theorem}\label{Cor_BinHarId1}
Let $m,n$ be positive integers with $m\geq n$. Then
\begin{multline*}
\sum_{k=0}^{n} \biggbin{m+k}{k} \biggbin{m}{k} \biggbin{n+k}{k} \biggbin{n}{k}
 \biggl[ 1+k \left(H_{m+k}^{(1)} +H_{m-k}^{(1)} + H_{n+k}^{(1)} 
 + H_{n-k}^{(1)} -4H_k^{(1)}\right) \biggr]\\[6pt]
\notag +\sum_{k=n+1}^{m} (-1)^{k-n} \biggbin{m+k}{k} \biggbin{m}{k} \biggbin{n+k}{k} \Big/ \biggbin{k-1}{n}
=(-1)^{m+n}.
\end{multline*}
\end{theorem}

\begin{theorem}\label{Cor_BinHarId2}
Let $l,m, n$ be positive integers with $l > m\geq n\geq\frac{l}{2}$ and $c_1, c_2 \in \mathbb{Q}$ some constants. Then
\begin{multline*}
\sum_{k=0}^{n} \biggbin{m+k}{k} \biggbin{m}{k} \biggbin{n+k}{k} \biggbin{n}{k} 
 \Biggl\{ \biggl[1+k \Bigl(H_{m+k}^{(1)} +H_{m-k}^{(1)} + H_{n+k}^{(1)} + H_{n-k}^{(1)} 
-4H_k^{(1)} \Bigr)\biggr]
\\[5pt]
 \cdot \biggl[c_1\left(H_{k+n}^{(1)} - H_{k+l-n-1}^{(1)}\right)  + c_2 \Bigl(H_{k+m}^{(1)} - H_{k+l-m-1}^{(1)}\Bigr)\biggr] 
-k\biggl[c_1\left(H_{k+n}^{(2)} - H_{k+l-n-1}^{(2)}\right)\\[5pt] 
+ c_2 \left(H_{k+m}^{(2)} - H_{k+l-m-1}^{(2)}\right)\biggr] \Biggr\}
+ \sum_{k=n+1}^{m} (-1)^{k-n} \biggbin{m+k}{k} \biggbin{m}{k} \biggbin{n+k}{k} \Big/ \biggbin{k-1}{n} 
\\[5pt]
\cdot \biggl[c_1\left(H_{k+n}^{(1)} - H_{k+l-n-1}^{(1)}\right) + c_2 \left(H_{k+m}^{(1)} - H_{k+l-m-1}^{(1)}\right)\biggr] = 0.
\end{multline*}
\end{theorem}
\noindent The remainder of this paper is spent proving Theorems \ref{Cor_BinHarId1} and \ref{Cor_BinHarId2}.

\section{Proofs}
We first develop two algebraic identities of which the binomial coefficient--harmonic sum identities are limiting cases.

\begin{theorem}\label{Thm_BinHarId1}
Let $x$ be an indeterminate and let $m,n$ positive integers with $m\geq n$. Then
\begin{multline}\label{BinHarIdx1}
\sum_{k=0}^{n} \biggbin{m+k}{k} \biggbin{m}{k} \biggbin{n+k}{k} \biggbin{n}{k}\\
\cdot \Biggl\{ \frac{-k}{(x+k)^2}
+  \frac{1+k \left(H_{m+k}^{(1)} +H_{m-k}^{(1)} + H_{n+k}^{(1)} + H_{n-k}^{(1)} -4H_k^{(1)}\right)}{x+k} \Biggr\} \\[6pt]
\raisetag{21pt} +\sum_{k=n+1}^{m} \frac{(-1)^{k-n}}{x+k} \biggbin{m+k}{k} \biggbin{m}{k} \biggbin{n+k}{k} \Big/ \biggbin{k-1}{n}
=\frac{x \ph{1-x}{n} \ph{1-x}{m}}{\ph{x}{n+1} \ph{x}{m+1}}.\phantom{=(2.1)}
\end{multline}
\end{theorem}
\vspace{6pt}
\begin{proof}
Using partial fraction decomposition we can write
\vspace{6pt}
\begin{equation*}
f(x):= \frac{x \ph{1-x}{n} \ph{1-x}{m}}{\ph{x}{n+1} \ph{x}{m+1}}= \frac{A}{x} + \sum_{k=1}^{n}
\biggl\{ \frac{B_k}{(x+k)^2} + \frac{C_k}{x+k} \biggr\} +  \sum_{k=n+1}^{m} \frac{D_k}{x+k}
\end{equation*}
for some $A, B_k, C_k$ and $D_k \in \mathbb{Q}$. We now isolate these coefficients by taking various limits of $f(x)$ as follows.
\begin{align*}
A = \lim_{x \to 0} x f(x) 
= \lim_{x \to 0}\frac{\ph{1-x}{n} \ph{1-x}{m}}{\ph{1+x}{n} \ph{1+x}{m}}
=1.
\end{align*}
\vspace{6pt}

\noindent For $1 \leq k \leq n$,
\begin{align*}
B_k = \lim_{x \to-k} (x+k)^2 f(x)
&=\lim_{x \to-k} \frac{x \ph{1-x}{n} \ph{1-x}{m}}{\ph{x}{k}^2 \ph{x+k+1}{n-k} \ph{x+k+1}{m-k}}\\[12pt]
&=\frac{-k \ph{k+1}{n} \ph{k+1}{m}}{\ph{-k}{k}^2 \ph{1}{n-k} \ph{1}{m-k}}\\[12pt]
&=\frac{-k \ph{k+1}{n} \ph{k+1}{m}}{(-1)^{2k} k!^2 (n-k)! (m-k)!}\\[12pt]
&=-k  \biggbin{m+k}{k} \biggbin{m}{k} \biggbin{n+k}{k} \biggbin{n}{k},\\[-3pt]
\end{align*}
and, using L'H\^{o}spital's rule,\\[-3pt]
\begin{align*}
C_k &=  \lim_{x \to-k} \frac{(x+k)^2 f(x) - B_k}{x+k}\\[12pt]
&=\lim_{x \to-k} \frac{d}{dx}\Biggl[(x+k)^2 f(x)\Biggr]\\[12pt]
&=\lim_{x \to-k}  \frac{d}{dx}\left[\frac{x \ph{1-x}{n} \ph{1-x}{m}}{\ph{x}{k}^2 \ph{x+k+1}{n-k} \ph{x+k+1}{m-k}}\right]\\[12pt]
&=\lim_{x \to-k} \left\{\left[\frac{ \ph{1-x}{n} \ph{1-x}{m}}{\ph{x}{k}^2 \ph{x+k+1}{n-k} \ph{x+k+1}{m-k}}\right]
\left[1-x \left(\sum_{s=1}^{n}(-x+s)^{-1}\right.\right.\right.\\[6pt]
&   \left.\left.\left.
\quad \phantom{=} +\sum_{s=1}^{m}(-x+s)^{-1} +\sum_{s=1}^{n-k}(x+k+s)^{-1}+\sum_{s=1}^{m-k}(x+k+s)^{-1}+2\sum_{s=0}^{k-1}(x+s)^{-1}\right)\right]\right\}\\[12pt]
&=\left[\frac{ \ph{1+k}{n} \ph{1+k}{m}}{\ph{-k}{k}^2 \ph{1}{n-k} \ph{1}{m-k}}\right]
\left[1+k \left(\sum_{s=1}^{n}(k+s)^{-1}+\sum_{s=1}^{m}(k+s)^{-1} +\sum_{s=1}^{n-k}(s)^{-1}\right.\right.\\
&\qquad \qquad \qquad\qquad \qquad\qquad \qquad\qquad \qquad \qquad \quad \;\:\:
\left.\left.+\sum_{s=1}^{m-k}(s)^{-1}+2\sum_{s=0}^{k-1}(-k+s)^{-1}\right)\right]\\[6pt]
&= \biggbin{m+k}{k} \biggbin{m}{k} \biggbin{n+k}{k} \biggbin{n}{k} \biggl[1+k \left(H_{m+k}^{(1)} +H_{m-k}^{(1)} + H_{n+k}^{(1)} + H_{n-k}^{(1)} -4H_k^{(1)}\right)\biggr].
\end{align*}
\vspace{9pt}

\noindent Similarly, for $n+1 \leq k \leq m$,
\begin{align*}
D_k = \lim_{x \to-k} (x+k) f(x)
&=\lim_{x \to-k} \frac{x \ph{1-x}{n} \ph{1-x}{m}}{\ph{x}{n+1} \ph{x}{k} \ph{x+k+1}{m-k}}\\[12pt]
&=\frac{-k \ph{k+1}{n} \ph{k+1}{m}}{\ph{-k}{n+1} \ph{-k}{k} \ph{1}{m-k}}\\[12pt]
&=(-1)^{k-n}  \biggbin{m+k}{k} \biggbin{m}{k} \biggbin{n+k}{k} \Big/ \biggbin{k-1}{n}.
\end{align*}

\end{proof}

\begin{theorem}\label{Thm_BinHarId2}
Let $x$ be an indeterminate and let $l,m, n$ be positive integers with $l > m\geq n\geq\frac{l}{2}$ and $c_1, c_2 \in \mathbb{Q}$ some constants. Then
\begin{multline}\label{BinHarIdx2}
\sum_{k=0}^{n} \frac{1}{x+k} \biggbin{m+k}{k} \biggbin{m}{k} \biggbin{n+k}{k} \biggbin{n}{k}
 \Biggl\{ \biggl[c_1\left(H_{k+n}^{(1)} - H_{k+l-n-1}^{(1)}\right) \\
+ c_2 \left(H_{k+m}^{(1)} - H_{k+l-m-1}^{(1)}\right)\biggr]
\cdot \biggl[ \frac{x}{x+k} + k \left(H_{m+k}^{(1)} +H_{m-k}^{(1)} + H_{n+k}^{(1)} + H_{n-k}^{(1)} -4H_k^{(1)}\right) \biggr]
\\[6pt]
\shoveleft{\qquad - k \biggl[c_1\left(H_{k+n}^{(2)} - H_{k+l-n-1}^{(2)}\right) 
 + c_2 \left(H_{k+m}^{(2)} - H_{k+l-m-1}^{(2)}\right)\biggr]\Biggr\}} \\
\shoveleft{+\sum_{k=n+1}^{m} \frac{(-1)^{k-n}}{x+k} \biggbin{m+k}{k} \biggbin{m}{k} \biggbin{n+k}{k} \Big/ \biggbin{k-1}{n}}\\[6pt]
\shoveleft  
\qquad \cdot \biggl[c_1\left(H_{k+n}^{(1)} - H_{k+l-n-1}^{(1)}\right)
+ c_2 \left(H_{k+m}^{(1)} - H_{k+l-m-1}^{(1)}\right)\biggr]
\\[6pt] 
=\frac{x \ph{1-x}{n} \ph{1-x}{m}}{\ph{x}{n+1} \ph{x}{m+1}}
\left[c_1 \sum_{s=l-n}^{n} (-x+s)^{-1} + c_2 \sum_{s=l-m}^{m} (-x+s)^{-1} \right].
\end{multline}
\end{theorem}
\vspace{6pt}
\begin{proof}
Using partial fraction decomposition we can write
\vspace{6pt}
\begin{align*}
f(x):&=\frac{x \ph{1-x}{n} \ph{1-x}{m}}{\ph{x}{n+1} \ph{x}{m+1}}
\left[c_1 \sum_{s=l-n}^{n} (-x+s)^{-1} + c_2 \sum_{s=l-m}^{m} (-x+s)^{-1} \right]\\[6pt]
&=\frac{A}{x} + \sum_{k=1}^{n} \biggl\{ \frac{B_k}{(x+k)^2} + \frac{C_k}{x+k} \biggr\} +  \sum_{k=n+1}^{m} \frac{D_k}{x+k}
\end{align*}
for some $A, B_k, C_k$ and $D_k \in \mathbb{Q}$. As in the proof of Theorem \ref{Thm_BinHarId1}, we isolate the coefficients $A, B_k, C_k$ and $D_k$ by taking various limits of $f(x)$.
For brevity, we first let
\vspace{6pt}
\begin{equation*}
T_{a}^{(r)}:= c_1 \sum_{s=l-n}^{n} (a+s)^{-r} + c_2 \sum_{s=l-m}^{m} (a+s)^{-r}
\end{equation*}
and
\vspace{6pt}
\begin{equation*}
U^{(r)} :=c_1\left(H_{k+n}^{(r)} - H_{k+l-n-1}^{(r)}\right) 
+ c_2 \left(H_{k+m}^{(r)} - H_{k+l-m-1}^{(r)}\right).
\end{equation*}

\noindent Then we have
\begin{align*}
A = \lim_{x \to 0} x f(x)
&= c_1 \lim_{x \to 0} \sum_{s=l-n}^{n} \frac{\ph{1-x}{n} \ph{1-x}{m}}{\ph{1+x}{n} \ph{1+x}{m} (s-x)}
+ c_2 \lim_{x \to 0} \sum_{s=l-m}^{m} \frac{\ph{1-x}{n} \ph{1-x}{m}}{\ph{1+x}{n} \ph{1+x}{m} (s-x)}\\[7pt]
&=c_1 \sum_{s=l-n}^{n} s^{-1} +c_2 \sum_{s=l-m}^{m} s^{-1}\\[7pt]
&=c_1\left(H_n^{(1)} - H_{l-n-1}^{(1)}\right)  + c_2 \left(H_m^{(1)} - H_{l-m-1}^{(1)}\right).
\end{align*}

\noindent For $1 \leq k \leq n$,
\begin{align*}
B_k = \lim_{x \to-k} (x+k)^2 f(x)
&=\lim_{x \to-k} \frac{x \ph{1-x}{n} \ph{1-x}{m}}{\ph{x}{k}^2 \ph{x+k+1}{n-k} \ph{x+k+1}{m-k}}
\;T_{-x}^{(1)}\\[7pt]
&=\frac{-k \ph{k+1}{n} \ph{k+1}{m}}{\ph{-k}{k}^2 \ph{1}{n-k} \ph{1}{m-k}}
\;T_{k}^{(1)}\\[7pt]
&=-k  \biggbin{m+k}{k} \biggbin{m}{k} \biggbin{n+k}{k} \biggbin{n}{k} U^{(1)}
\end{align*}
and
\begin{align*}
C_k &=\lim_{x \to-k} \frac{d}{dx}\Biggl[(x+k)^2 f(x)\Biggr]\\[7pt]
&=\lim_{x \to-k}  \frac{d}{dx}\left[\frac{x \ph{1-x}{n} \ph{1-x}{m}}{\ph{x}{k}^2 \ph{x+k+1}{n-k} \ph{x+k+1}{m-k}}
\;T_{-x}^{(1)} \right] \\[7pt]
&=\lim_{x \to-k} \Biggl\{ \Biggl[\frac{ \ph{1-x}{n} \ph{1-x}{m}}{\ph{x}{k}^2 \ph{x+k+1}{n-k} \ph{x+k+1}{m-k}}\Biggr]
\Biggl[\;x\; T_{-x}^{(2)} +\;T_{-x}^{(1)} -x\;T_{-x}^{(1)} \\[4pt]
& \phantom{=} \qquad \qquad 
\cdot \Biggl(\sum_{s=1}^{n}(-x+s)^{-1}
 +\sum_{s=1}^{m}(-x+s)^{-1}
+\sum_{s=1}^{n-k}(x+k+s)^{-1}
 \\[4pt] &\qquad \qquad  \qquad \qquad \qquad \qquad \qquad \qquad
+\sum_{s=1}^{m-k}(x+k+s)^{-1}+2\sum_{s=0}^{k-1}(x+s)^{-1}\Biggr)\Biggr] \Biggr\}\\[7pt]
&=\left[\frac{ \ph{1+k}{n} \ph{1+k}{m}}{\ph{-k}{k}^2 \ph{1}{n-k} \ph{1}{m-k}}\right]
\Biggl[-k T_{k}^{(2)}+ T_{k}^{(1)}
\left(1+k \left(\sum_{s=1}^{n}(k+s)^{-1}+\sum_{s=1}^{m}(k+s)^{-1} \right.\right.\\[4pt]
&\qquad \qquad \qquad \qquad \qquad \qquad \qquad  \left.\left.
+\sum_{s=1}^{n-k}(s)^{-1} +\sum_{s=1}^{m-k}(s)^{-1}+2\sum_{s=0}^{k-1}(-k+s)^{-1}\right)\right)\Biggr]\\[7pt]
&= \biggbin{m+k}{k} \biggbin{m}{k} \biggbin{n+k}{k} \biggbin{n}{k} \\*[4pt]
&  \qquad \qquad
 \cdot \biggl[-k \,U^{(2)} +\bigg(1+k \biggl(H_{m+k}^{(1)} +H_{m-k}^{(1)} + H_{n+k}^{(1)} 
+ H_{n-k}^{(1)} -4H_k^{(1)}\biggr)\biggr) U^{(1)} \biggr].
\end{align*}
\noindent For $n+1 \leq k \leq m$,
\begin{align*}
D_k = \lim_{x \to-k} (x+k) f(x)
&=\lim_{x \to-k} \frac{x \ph{1-x}{n} \ph{1-x}{m}}{\ph{x}{n+1} \ph{x}{k} \ph{x+k+1}{m-k}}
\;T_{-x}^{(1)}\\[12pt]
&=\frac{-k \ph{k+1}{n} \ph{k+1}{m}}{\ph{-k}{n+1} \ph{-k}{k} \ph{1}{m-k}}
\;T_{k}^{(1)}\\[12pt]
&=(-1)^{k-n} \, U^{(1)} \, \biggbin{m+k}{k} \biggbin{m}{k} \biggbin{n+k}{k} \Big/ \biggbin{k-1}{n}.
\end{align*}
\end{proof}

\begin{proof}[Proofs of Theorems \ref{Cor_BinHarId1} and \ref{Cor_BinHarId2}]
Multiply both sides of (\ref{BinHarIdx1}) and (\ref{BinHarIdx2}) respectively by $x$ and take the limit as $x \to \infty$.
\end{proof}


\begin{thebibliography}{999}

\bibitem{AEOZ} S. Ahlgren, S. B. Ekhad, K. Ono, D. Zeilberger, \emph{A binomial coefficient identity associated to a conjecture of Beukers}, Electron. J. Combin. \textbf{5} (1998), Research Paper 10, 1.

\bibitem{AO} S. Ahlgren, K. Ono, \emph{A Gaussian hypergeometric series evaluation and Ap{\'e}ry number congruences}, J. Reine Angew. Math. \textbf{518} (2000), 187--212.

\bibitem{C} W. Chu, \emph{A binomial coefficient identity associated with Beukers' conjecture on Ap\'ery numbers}, Electron. J. Combin. \textbf{11} (2004), no. 1, Note 15, 3.

\bibitem{DMC} D. McCarthy, \emph{$p$-adic hypergeometric series and supercongruences}, Ph.D thesis, University College Dublin, 2010.

\bibitem{DMC1} D. McCarthy, \emph{Extending Gaussian hypergeometric series to the $p$-adic setting}, Int. J. Number Theory, accepted for publication, 24 pages.

\bibitem{DMC2} D. McCarthy, \emph{On a supercongruence conjecture of Rodriguez-Villegas}, Proc. Amer. Math. Soc. \textbf{140} (2012), 2241--2254.

\end{thebibliography}
\end{document}